\documentclass[reqno,twoside,a4paper]{amsart}

\usepackage{amscd}
\usepackage{amssymb}
\usepackage{amsxtra}
\usepackage{longtable}
\usepackage{calc}
\usepackage{epsfig}

\allowdisplaybreaks[2]

\theoremstyle{plain}
\newtheorem{thm}{Theorem}
\newtheorem{cor}{Corollary}
\newtheorem{prop}{Proposition}
\newtheorem{lem}{Lemma}

\theoremstyle{definition}
\newtheorem{defn}{Definition}

\theoremstyle{remark}
\newtheorem{rmk}{Remark}
\newtheorem{rmks}{Remarks}
\newtheorem{exmp}{Example}
\newtheorem{exmps}{Examples}

\renewcommand{\epsilon}{\varepsilon}
\renewcommand{\kappa}{\varkappa}
\renewcommand{\theta}{\vartheta}

\newcommand{\IC}{\ensuremath{\mathbb C}}

\newcommand{\IH}{\ensuremath{\mathbb H}}
\newcommand{\IN}{\ensuremath{\mathbb N}}
\newcommand{\IQ}{\ensuremath{\mathbb Q}}

\newcommand{\IZ}{\ensuremath{\mathbb Z}}

\newcommand{\ch}{\ensuremath{\mathrm{ch}}}
\newcommand{\td}{\ensuremath{\mathrm{td}}}

\newcommand{\ellip}{\ensuremath{\mathrm{Ell}}}

\newcommand{\Sym}{\ensuremath{\mathrm{S}}}

\newcommand{\varcard}[1]{\ensuremath{\operatorname{card}}}
\newcommand{\ord}[1]{\ensuremath{\operatorname{ord}}}

\newcommand{\HH}{\ensuremath{\mathrm H}}
\newcommand{\KK}{\ensuremath{\mathrm K}}

\newcommand{\diff}{\ensuremath{\mathrm d}}

\newcommand{\Tang}{\ensuremath{\mathcal T}}
\newcommand{\TangX}{\ensuremath{\Tang_X}}

\newcommand{\pt}{\ensuremath{\mathrm{point}}}

\def\Z{\IZ}
\def\Q{\IQ}

\title{On the Elliptic Genus of Generalised Kummer varieties}
\author[M.\ A.\ Nieper-Wi\ss kirchen]{Marc A.~Nieper-Wi\ss kirchen}
\address{Mathematisches Institut der Univ.\ zu K\"oln \\
  Weyertal 86--90 \\ 50931 K\"oln \\ Germany}
\email{mnieper@mi.uni-koeln.de}
\date{\today}

\begin{document}

\begin{abstract}
  Borisov and Libgober (\cite{borisov02}) recently proved a conjecture of
  Dijkgraaf,
  Moore, Verlinde, and Verlinde (see~\cite{dijkgraaf97}) on the elliptic genus
  of a Hilbert scheme of points on a surface. We show how their result can be
  used together with our work on complex genera of generalised Kummer
  varieties~\cite{kummer} to deduce the following formula on the elliptic
  genus of a generalised Kummer variety $A^{[[n]]}$ of dimension $2(n - 1)$:
  \begin{gather*}
    \ellip(A^{[[n]]}) = n^4 \psi^{-2} \cdot \left(\left.\psi^2\right|_{-2}
    V(n)\right).
  \end{gather*}
  Here $\psi(z, \tau) := 2 \pi i \frac{\theta(-z, \tau)}{\theta'(0,
  \tau)}$ is the weak Jacobi form of weight $-1$ and index $\frac 1 2$, and
  $V(n)$ is the Hecke operator sending Jacobi forms of index $r$ to Jacobi
  forms of index $nr$ (see~\cite{ez85}).
\end{abstract}

\maketitle

\section{Introduction}

The elliptic genus can be defined for every compact complex manifold as a
function of a complex variable $\tau \in \IH$ and a complex variable $z \in
\IC$.  For a Calabi-Yau manifold (Calabi-Yau in the (very weak) sense that its
first Chern class vanishes up to torsion) the elliptic genus turns out to be a
weak Jacobi form of weight zero and index $d/2$, where $d$ is the dimension of
the manifold. A proof of this fact can be found in~\cite{gritsenko99}, and
also in~\cite{borisov00}.

The elliptic genus is a generalization of Hirzebruch's $\chi_y$-genus. It was
shown in~\cite{borisov00} that in fact the elliptic genus of a Calabi-Yau
generally encodes more information than the $\chi_y$-genus provided the
dimension of the manifold is twelve or greater than thirteen.

Compact hyperk\"ahler (see for example~\cite{huybrechts99}) manifolds are in
particular Calabi-Yau manifolds in the sense given above. There are two main
series of examples of compact hyperk\"ahler manifolds, the Hilbert schemes of
points on surfaces and the generalised Kummer varieties
(see~\cite{beauville83}). The $\chi_y$-genera for these manifolds have been
known for quite some time (see~\cite{goettsche93}) since they
depend only on the Hodge numbers of these manifolds. One is naturally led
to the question of generalising these results to elliptic genera.

Dijkgraaf, Moore, Verlinde and Verlinde showed in~\cite{dijkgraaf97} that the
elliptic genus of a symmetric product of a manifold can be expressed in terms
of the elliptic genus of the manifold itself. They conjectured that their
formula also holds for the Hilbert scheme of points on surfaces since these
Hilbert schemes can be seen as resolutions of symmetric powers of surfaces.
Very recently, Borisov and Libgober managed to prove that conjecture
(see~\cite{borisov02}) with methods of toric geometry.

In~\cite{kummer} we expressed the value of a genus on a generalised Kummer
variety in terms of twisted versions of this genus evaluated on the  Hilbert
schemes of points on a surface $X$ with $c_1(X)^2 \neq 0$. In this work it is
shown that the result of~\cite{borisov02} on the elliptic genera of these
Hilbert schemes can be used to calculate their twisted elliptic genera. We do
this and are led to a formula for the elliptic genus of a generalised Kummer
variety. By specialising our result we reproduce G\"ottsche's and Soergel's
formula for the $\chi_y$-genus of the generalised Kummer varieties.

One last word on the notation: The symbol $\sum_{a|n}$ stands for summing over
the~\emph{positive} divisors of the integer $n$.

\section{Complex genera}

\subsection{Some notes on complex genera and the elliptic genus}
Let $\Omega := \Omega^{\mathrm{U}} \otimes \IQ$ denote the (rational) complex
cobordism ring. By a result of Milnor (\cite{milnor60}), it is generated by
the cobordism classes $[X]$ of all complex manifolds $X$, and two complex
manifolds $X$ and $Y$ lie in the same cobordism class if and only if they have
the same Chern numbers, i.e.\ the Chern numbers determine the cobordism class
and vice versa. Recall that the sum in the ring is induced by the disjoint
union of manifolds, and the product by the cartesian product of manifolds.

By a \emph{complex genus $\phi$ with values in a $\Q$-algebra $R$} we want to
understand a $\Q$-algebra homomorphism $\phi: \Omega \to R$. Examples of
complex genera can be constructed in the following way:
\begin{exmp}
  \label{ex:phigenus}
  Let $\mathfrak C$ be the category of complex manifolds and $\mathfrak A$ the
  category of abelian groups.  Let $R$ be any $\Q$-algebra and $\Phi$ a natural
  transformation from the functor $\mathfrak C \to \mathfrak A, X \mapsto
  \KK^0(X) \otimes \Q$ to the functor $\mathfrak C \to \mathfrak A, X \mapsto
  (\KK^0(X) \otimes R)^{\times}$ that associates to each $X$ in the category
  of compact complex manifolds a group homomorphism $\Phi_X: \KK^0(X) \to
  (\KK^0(X) \otimes R)^{\times}$. For every such natural transformation,
  \begin{gather}
    \phi: X \mapsto \chi(X, \Phi_X(\TangX))
  \end{gather}
  defines a complex genus.
\end{exmp}

\begin{proof}
  It follows directly from the naturality of $\Phi_X$ that $\phi(X + Y) =
  \phi(X) + \phi(Y)$. Obviously, $\phi(\emptyset) = 0$, where ``$+$'' denotes
  the disjoint union of manifolds. Furthermore,
  \begin{align*}
    \phi(X \times Y)
    = & \chi(X \times Y, \Phi_{X \times Y}(\Tang_{X \times Y}))
    \\
    = & \chi(X \times Y, \Phi_{X \times Y}(\Tang_X \boxplus \Tang_Y))
    \\
    = & \chi(X \times Y, \Phi_{X \times Y}(\Tang_X) \otimes
    \Phi_{X \times Y}(\Tang_Y))
    \\
    = & \chi(X \times Y, \Phi_X(\Tang_X) \boxtimes \Phi_Y(\Tang_Y))
    \\
    = & \chi(X, \Phi_X(\Tang_X)) \cdot \chi(Y, \Phi_Y(\Tang_Y))
    \\
    = & \phi(X) \cdot \phi(Y)
    \\
    \intertext{and}
    \phi(\pt) = & \chi(\pt, \Phi_\pt(0)) = \chi(\pt, 1) = 1.
  \end{align*}
 
  Finally, we have to show that $\phi$ depends only on the characteristic
  numbers of $X$. By the Hirzebruch-Riemann-Roch formula, we have $\phi(X) =
  \int_X \ch(\Phi_X(\TangX)) \td(\TangX)$, i.e.\ it suffices to show that
  $\ch(\Phi_X(E))$ depends only on the Chern classes of $E$ for any bundle $E$
  on $X$. This can be seen as follows: Obviously, $E \mapsto \ch(\Phi_X(E))$
  is a natural transformation from the functor $X \mapsto \KK^0(X) \otimes \Q$
  to the functor $X \mapsto \HH^0(X) \otimes R$. By the universality property
  of the Chern classes it follows that $\ch(\Phi_X(E))$ can be given as an
  universal polynomial (depending only on $\Phi_X$) in the Chern classes of
  $E$ over $R$ (see for example~\cite{bott82})
\end{proof}

The definition of the \emph{elliptic genus} is an application to the concept
outlined in the example. From now on let use write $q := e^{2 \pi i \tau}$ for
$\tau \in \IH$, and $y^{\frac 1 2} := e^{\pi i z}$ for $z \in \IC$, i.e.\ $y =
e^{2 \pi i z}$.
\begin{defn}
  \label{def:ebundle}
  For every complex manifold $X$ of dimension $d$, we define the following
  virtual bundle
  \begin{gather}
    \mathcal E_X(z, \tau) := \Phi_X(\TangX, z, \tau)
    \in \KK^0(X)[[q]][y^{\pm \frac 1 2}]
    \\
    \intertext{on $X$, where}
    \begin{aligned}
      \Phi_X(E, z, \tau)
      & = y^{-d/2} \bigotimes_{n = 1}^\infty
      \left(\bigwedge\nolimits_{-y q^{n - 1}} E^*
        \otimes \bigwedge\nolimits_{-y^{-1}q^n} E
        \otimes \Sym_{q^n} E^* \otimes \Sym_{q^n} E\right)
      \\
      & \in \KK^0(X)[[q]][y^{\pm \frac 1 2}],
    \end{aligned}
  \end{gather}
  for any bundle $E$ on $X$. Here, $\bigwedge_t(E) := \bigoplus_{i
    = 0}^\infty \bigwedge^i(E) t^i$ and $\Sym_t(E) := \bigoplus_{i = 0}^\infty
  \Sym^i(E) t^i$ for any bundle $E$ denote the generating functions for the
  exterior respective symmetric powers of $E$.
\end{defn}

\begin{defn}[see~\cite{borisov00}]
  Let $X$ be a compact complex manifold of dimension $d$. The \emph{elliptic
    genus $\ellip(X, z, \tau)$ of $X$} is defined to be the Euler
  characteristic $\chi(X, \mathcal E_X(z, \tau)) \in \Q[[q]][y^{\pm \frac 1
    2}]$ of the bundle $\mathcal E_X(z, \tau)$ on $X$.
\end{defn}

The natural transformation $X \mapsto \Phi_X$ of Definition~\ref{def:ebundle}
is in fact a natural
transformation as in Example~\ref{ex:phigenus}. Therefore, $X \mapsto
\ellip(X)$ induces in fact a complex genus.

\begin{rmk}
  The elliptic genus is a generalization of the $\chi_y$-genus of
  Hirzebruch. More precisely,
  \begin{gather}
    \chi_{-y}(X) = y^{d/2} \lim_{\tau \to i \infty} \ellip(X, z, \tau)
    = y^{d/2} \left.\ellip(X, z, \tau)\right|_{q = 0} \in \Q[y],
  \end{gather}
  which follows directly from the definitions.
\end{rmk}

By Hirzebruch's theory of genera and multiplicative sequences
(\cite{hirzebruch56}), the $R$-valued complex genera $\phi$ are in one-to-one
correspondence with the formal power series $f_\phi \in R[[x]]$ over $R$ with
constant coefficient one. The correspondence is given as follows:
\begin{equation}
  \phi(X) = \int_X \prod_{i = 1}^d f_\phi(x_i)
\end{equation}
for all complex manifolds $X$, where $d$ is the dimension of $X$, and
$x_1, \dots, x_d$ are the Chern roots of its tangent bundle.

The representation of the elliptic genus in terms of power series is as
follows:
\begin{prop}
  Let $x_1, \dots, x_d$ be the Chern roots of $X$. We have
  \begin{gather}
    \ellip(X, z, \tau) = \int_X \prod_{i = 1}^d \left(
      x_i \frac{\theta\left(\frac{x_i}{2 \pi i} - z, \tau\right)}
      {\theta\left({\frac{x_i}{2 \pi i}}\right)}\right),
  \end{gather}
  where Jacobi's theta function is defined as
  \begin{gather}
    \theta: \IC \times \IH \to \IC,
    (z, \tau) \mapsto 2 q^{1/8} \sin(\pi z) \prod_{n = 1}^\infty
    \left((1-q^n)(1-q^n y)(1-q^n y^{-1})\right).
  \end{gather}
\end{prop}

\begin{proof}
  See~\cite{borisov00}.
\end{proof}

However, the power series in $x_i$ on the right hand side of the product sign
has not constant coefficient one. Nevertheless, we can achieve this by
following an idea proposed in~\cite{hbj92}, where it has been worked out on
the $\chi_y$-genus: Let us define
\begin{gather}
  \psi(z, \tau) := \lim_{x \to 0} x \frac {\theta(\frac x{2\pi i} - z, \tau)}
  {\theta(\frac x{2 \pi i}, \tau)}
  = 2 \pi i \frac{\theta(-z, \tau)}{\theta'(0, \tau)} \in
  \Q[[q]][y^{\pm \frac 1 2}],
\end{gather}
where $\theta'(0, \tau) :=
\left.\frac{\partial}{\partial z}\right|_{z = 0} \theta(z, \tau)$.
Furthermore, we set
\begin{gather}
  f(z, \tau)(x) := \frac x \psi \frac{\theta(\frac{x \psi}{2 \pi i} - z,
    \tau)} {\theta(\frac{x \psi}{2 \pi i}, \tau)}
  \in \Q[[q]][y^{\pm \frac 1 2}][[x]].
\end{gather}
This is a power series in $x$ with constant coefficient one.
With this definition, we have
\begin{gather}
  \label{equ:ellipf}
  \ellip(X, z, \tau) = \int_X \prod_{i = 1}^d f(z, \tau)(x_i),
\end{gather}
which can be seen by expanding $f(z, \tau)$ as a Taylor series in $x$ and
looking at the degree $d$ term.

\begin{rmk}
  The power series $\psi$ has the following representation as an infinite
  pro\-duct:
  \begin{gather}
    \psi(z, \tau) = \left(y^{- \frac 1 2} - y^{\frac 1 2}\right)
    \prod_{n = 1}^\infty \frac{(1 - q^n y)(1 - q^n y^{-1})}{(1 - q^n)^2}.
  \end{gather}
  This can be seen for example by considering the product representations of
  $\theta(\tau, z)$ and $\theta'(0, \tau)$ given in~\cite{chandra85}. 
\end{rmk}

\subsection{Twisted genera and the twisted elliptic genus}

Following our work on complex genera of generalised Kummer varieties
(\cite{kummer}), we define \emph{twisted genera of compact complex manifolds}
as follows:
\begin{defn}
  Let $\phi$ be a complex genus with values in the $\Q$-algebra $R$. By
  $\phi_t$ we denote the genus with values in $R[t]$ that is given by
  \begin{equation}
    \phi_t(X) := \int_X \prod_{i = 1}^d \left(f_\phi(x_i)
    e^{t x_i}\right)
  \end{equation}
  for any complex manifold $X$ of dimension $d$, and with Chern roots $x_1,
  \dots, x_d$.
\end{defn}

\begin{rmk}
  Obviously, if a the genus $\phi$ is given by a natural transformation $\Phi$
  as in Example~\ref{ex:phigenus}, we have
  \begin{gather}
    \phi_t(X) = \chi(X, \Phi_X(\Tang_X) \otimes K_X^{-t}).
  \end{gather}
  for any integer $t$. Of course, this can also be written as
  \begin{gather}
    \phi_t(X) = \chi(X, \Phi_{t, X}(\TangX))
  \end{gather}
  with the natural transformation $\Phi_t$ defined by $\Phi_{t, X}(E) :=
  \Phi_X(E) \otimes (\det E)^t$ for any bundle $E$ on $X$.
\end{rmk}

In general, the twisted version of a complex genus holds more information than
the untwisted genus itself. In the case of the elliptic genus, however, this
is not true, i.e.~the elliptic genus encodes all the information of its
twisted version.

As we are interested only in manifolds of even dimension for the purpose of
proving our main theorem, let $X$ be of even dimension $d = 2 n$ from now on.
\begin{prop}
  \label{prop:twistell}
  The elliptic genus twisted by $t \psi$, t an integer, of the compact complex
  manifold $X$ can be expressed by the usual elliptic genus via
  \begin{gather}
    \ellip_{t \psi}(X, z, \tau) = y^{2 n t} q^{n t^2} \ellip(X, z + t \tau,
    \tau).
  \end{gather}
\end{prop}

\begin{proof}
  We will make use of the following transformation property of the theta
  function (see for example~\cite{chandra85}):
  \begin{gather}
    \theta(z + \tau, \tau) = - q^{- \frac 1 2} y^{-1} \theta(z, \tau)
  \end{gather}
  where $q^{- \frac 1 2}$ stands for $\exp(- \pi i \tau)$.
  It follows by definition that
  \begin{gather}
    f(z, \tau, x) = - q^{\frac 1 2} y e^{-x \psi} f(x, z + \tau, \tau).
  \end{gather}
  Now we use~\eqref{equ:ellipf} and find that
  \begin{gather}
    \ellip_{t \psi}(X, z, \tau) = q^n y^{2n} \ellip_{(t - 1) \psi}(X, z +
    \tau, \tau)
  \end{gather}
  by the definition of twisted genera. The proposition then follows by
  induction on $t$.
\end{proof}

From this we deduce the following proposition on the Fourier coefficients of
twisted elliptic genera:
\begin{prop}
  Let $t$ be an integer and $c_t(m, l)$ for $m, 2 l \in \IZ$ be such that
  \begin{gather}
    \sum_{m, 2 l \in \Z} c_t(m, l) q^m y^l = \ellip_{t \psi}(X, z, \tau).
  \end{gather}
  These coefficients are connected to the coefficients $c(m, l) := c_0(m, l)$
  of the untwisted elliptic genus via
  \begin{gather}
    c_t(m, l) = c(m - l t + n t^2, l - 2 t).
  \end{gather}
\end{prop}

\begin{proof}
  By Proposition~\ref{prop:twistell}, we have
  \begin{multline*}
    \sum_{m, 2l \in \Z} c_t(m, l) q^m y^l = \ellip_{t \psi}(X, z, \tau)
    = y^{2nt} q^{n t^2} \ellip(X, z + t \tau, \tau)
    \\
    = y^{2nt} q^{n t^2} \sum_{m, 2l \in Z} c(m, l) q^m (yq^t)^{l}
    = \sum_{m, 2l \in \Z} c(m, l) q^{m + l t + n t^2} y^{l + 2nt}
    \\
    = \sum_{m, 2l \in \Z} c(m, l - 2nt) q^{m + l t - nt^2} y^l
    = \sum_{m, 2l \in \Z} c(m - lt + n t^2, l - 2nt) q^m y^l.    
  \end{multline*}
  Now compare coefficients.
\end{proof}

\section{The generalised Kummer varieties}

\subsection{Definition of generalised Kummer varieties}
Let $X$ be a smooth projective surface over the field of complex numbers. For
every nonnegative integer we denote by $X^{[n]}$ the Hilbert scheme of
zero-dimensional subschemes of $X$ of length $n$. By a result of
Fogarty~\cite{fogarty68}, this scheme is smooth and projective of dimension
$2n$. It can be viewed as a resolution $\rho: X^{[n]} \to X^{(n)}$ of
the $n$-fold symmetric product $X^{(n)} := X^n/\mathfrak{S}_n$ of $X$.
The morphism
$\rho$, sending closed points, i.e.\ subschemes of $X$, to their support
counting multiplicities, is called the Hilbert-Chow morphism.

Let us briefly recall the construction of the generalised Kummer
varieties introduced by Beauville~\cite{beauville83}.
Let $A$ be an abelian surface and $n > 0$. There is an obvious summation
morphism $A^{(n)} \to A$. We denote its composition with the Hilbert-Chow
morphism $\rho: A^{[n]} \to A^{(n)}$ by $\sigma: A^{[n]} \to A$.
\begin{defn}
 The \emph{$n^{\textrm{th}}$ generalised Kummer variety $A^{[[n]]}$}
 is the fibre of $\sigma$ over $0 \in A$.
\end{defn}

\begin{rmk}
  It was Beauville who showed among other things in~\cite{beauville83} that
  the $n^{\mathrm{th}}$ generalised Kummer variety is smooth projective, and
  irreducible holomorphic symplectic of dimension $2(n - 1)$.
\end{rmk}

\subsection{Complex genera of generalised Kummer varieties}

Let
\begin{gather}
  K(p) := \sum_{n = 1}^\infty \left[A^{[[n]]}\right] \cdot p^n \in \Omega[[p]]
\end{gather}
be the generating series of the complex cobordism classes of the generalised
Kummer varieties and similarly
\begin{gather}
  H_X(p) := \sum_{n = 0}^\infty \left[X^{[n]}\right] \cdot p^n \in \Omega[[p]]
\end{gather}
the generating series to the Hilbert schemes of points on a surface $X$.

In~\cite{kummer} we proved the following theorem:
\begin{thm}
  \label{thm:kummergen}
  Let $\phi: \Omega \to R$ be a complex genus with values in the $\Q$-algebra
  $R$. For every smooth projective surface $X$ with $\int_X c_1(X)^2 \neq
  0$,
  \begin{equation}
    \phi(K(p)) = \frac 1{c_1(X)^2} \left(p \frac{\diff}{\diff p}\right)^2
    \ln \frac{\phi_1(H_X(p)) \phi_{-1}(H_X(p))}{\phi(H_X(p))^2}.
  \end{equation}
\end{thm}

\begin{proof}
  See~\cite{kummer}.
\end{proof}

Taking a look at the proofs of~\cite{kummer} or at the results
of~\cite{lehn01}, we see that $\ln \phi_t(H_X(p))$ is a polynomial in $t$ of
degree $2$, and therefore
\begin{gather}
  \ln \frac{\phi_1(H_X(p)) \phi_{-1}(H_X(p))}{\phi(H_X(p))^2}
  = \left.\frac{\diff^2}{\diff t^2}\right|_{t = 0} \ln \phi_t(H_X(p)).
\end{gather}
So, the genus $\phi$ of $K(p)$ can be expressed as
\begin{gather}
  \label{equ:aphigen}
  \phi(K(p)) = 2 \left(p \frac{\diff}{\diff p}\right)^2 A_\phi(p)
\end{gather}
for any smooth projective surface $X$ with $c_1(X)^2 \neq 0$. Here,
\begin{gather}
  A_\phi(p) := \frac 1{2 c_1(X)^2} \left.\frac{\diff^2}{\diff t^2}\right|_{t =
    0} \ln \phi_t(H_X(p)),
\end{gather}
which is independent of $X$. We shall use the formula~\eqref{equ:aphigen} to
calculate the elliptic genera of the generalised Kummer varieties.

\section{Elliptic genera of generalised Kummer varieties}

In the following proposition, we shall calculate twisted elliptic genera of
Hilbert schemes of points on surfaces. We shall make use of the following
recent theorem of Borisov and Libgober (\cite{borisov02}):
\begin{thm}
  Let $X$ be a complex projective surface. Let $c(m, l)$ for $m, 2 l \in \Z$
  be the Fourier coefficients of the elliptic genus of $X$, i.e.~$\ellip(X, z,
  \tau) = \sum_{m, 2l \in \IZ} c(m, l) q^m y^l$. Then the elliptic genera of
  the Hilbert schemes $X^{[n]}$ of zero-dimensional subschemes of length $n$ on
  $X$ are given by:
  \begin{gather}
    \ellip(H_X(p), z, \tau) := \sum_{n = 0}^\infty p^n \ellip(X^{[n]}, z,
    \tau) = \prod_{i = 1}^\infty \prod_{m, 2l \in \IZ} (1 - p^i q^m y^l)^{-
    c(m i, l)}.
  \end{gather}
\end{thm}
 
\begin{proof}
  See~\cite{borisov02}.
\end{proof}

\begin{prop}
  \label{prop:elliptpsi}
  As a corollary to the previous theorem, we have
  \begin{gather}
    \ellip_{t \psi} (H_X(p), z, \tau) := \sum_{n = 0}^\infty p^n
    \ellip_{t \psi} (X^{[n]}, z, \tau) = \prod_{i = 1}^\infty \prod_{m, 2l \in
      \Z} (1 - p^i y^l q^m)^{- c_{it}(mi, l)},
    \intertext{for every integer $t$, i.e.}
    \ln \ellip_{t \psi} (H_X(p), z, \tau) = \sum_{i, k = 1}^\infty
    \sum_{m, 2l \in \IZ} \frac 1 k c_{it} (mi, l) p^{ik} y^{lk} q^{mk}.
  \end{gather}
\end{prop}

\begin{proof}
  This proof consists just of a straight-forward calculation:
  \begin{align*}
    \phantom = & \sum_{n = 0}^\infty p^n \ellip_{t \psi}(X^{[n]}, z, \tau) =
    \sum_{n = 0}^\infty p^n y^{2nt} q^{nt^2} \ellip(X^{[n]}, z + t \tau, \tau)
    \\
    = & \prod_{i = 1}^\infty \prod_{m, 2l \in \IZ} (1 - (py^{2t}q^{t^2})^i q^m
    (yq^t)^l)^{- c(mi, l)}
    \\
    = & \prod_{i = 1}^\infty \prod_{m, 2l \in \IZ} (1 - p^i q^{m + lt + it^2}
    y^{l + 2it})^{- c(mi, l)}
    \\
    = & \prod_{i = 1}^\infty \prod_{m, 2l \in \IZ} (1 - p^i y^l q^{m + l t - i
      t^2})^{- c(mi, l - 2it)}
    \\
    = & \prod_{i = 1}^\infty \prod_{m, 2l \in \IZ} (1 - p^i y^l q^m)^{- c(mi -
      lit + i^2 t^2, l - 2it)}
    \\
    = & \prod_{i = 1}^\infty \prod_{m, 2l \in \IZ} (1 - p^i y^l q^m)^{-
      c_{it}(mi, l)}.
  \end{align*}
\end{proof}

\begin{lem}
  \label{lem:csurft}
  Let $X$ be a smooth projective surface. Let $c_t(m, l)$ for $m, 2 l \in \Z$
  be the Fourier coefficients of the twisted elliptic genus $\ellip_{t \psi}$
  of $X$. We have
  \begin{gather}
    \left. \frac{\diff^2}{\diff t^2}\right|_{t = 0} c_t(m, l) = u(m,
    l) \cdot c_1(X)^2,
  \end{gather}
  where $u(m, l) \in \IQ$ for $m, 2l \in \IZ$ is such that
  \begin{gather}
    \sum_{m, 2l \in \IZ} u(m, l) q^m y^l = \psi^2(z, \tau).
  \end{gather}
\end{lem}

\begin{proof}
  Such a statement is actually true for any genus on a surface $X$:
  If $\phi$ is any complex genus, then
  \begin{gather*}
    \tag{$\ast$}
    \left. \frac{\diff^2}{\diff t^2}\right|_{t = 0} \phi_t(X) = c_1(X)^2,
  \end{gather*}
  which follows from the fact that $\phi_t(X)$ is polynomial of degree two in
  $t$ with leading coefficient $c_1(X)^2/2$.
  From ($\ast$) the lemma follows by the chain rule.
\end{proof}

\begin{prop}
  For every smooth projective surface $X$ and $u(m, l)$ for $m, 2l \in \IZ$
  as in Lemma~\ref{lem:csurft} we have:
  \begin{gather}
    \left.\frac{\diff^2}{\diff t^2}\right|_{t = 0} \ln \ellip_{t \psi}(H_X(p),
    z, \tau) = c_1(X)^2 \sum_{i, k = 1}^\infty \sum_{m, 2l \in \IZ} i^2 \frac
    1 k p^{ik} q^{mk} y^{lk} u(mi, l).
  \end{gather}
\end{prop}

\begin{proof}
  This proposition follows immediately from Proposition~\ref{prop:elliptpsi}
  and Lemma~\ref{lem:csurft}.
\end{proof}

Now we can use that proposition together with the theorem on complex
genera of generalised Kummer varieties to prove our main theorem:

\begin{thm}
  Let $J_{k, r}$ for $k, 2r \in \IZ$ be the space of weak Jacobi forms of
  weight $k$ and index $r$. Following~\cite{ez85} we define for every $i \in
  \IN$ a Hecke operator
  \begin{gather}
    V(i): J_{k, r} \to J_{k, ir}, \phi \mapsto \left.\phi\right|_k V(i)
  \end{gather}
  by
  \begin{gather}
    (\left.\phi\right|_k V(i))(z, \tau) := \sum_{a|i} a^{k - 1} \sum_{m, 2l \in
    \IZ} c\left(\frac{mi}a, l\right) q^{am} y^{al},
  \end{gather}
  where the $c(m, l)$ are the Fourier coefficients of $\phi$, i.e.~$\phi(z,
  \tau) = \sum_{m, 2l \in \Z} c(m, l) q^m y^l$.
 
  The following statements hold true: $\psi^2$ is a weak Jacobi form of weight
  $-2$ and index $1$, $\psi^{-2} \left(\left.\psi^2\right|_{-2} V(i)\right)$
  is a weak Jacobi form of weight $0$ and index $i - 1$, and the elliptic
  genera of the generalised Kummer varieties are given by
  \begin{gather}
    \ellip(K(p)) = \sum_{i = 1}^\infty \ellip(A^{[[n]]}) p^n
    = \psi^{-2} \sum_{i = 1}^\infty i^4
    \left(\left.\psi^2\right|_{-2} V(i)\right) p^i.
  \end{gather}
\end{thm}

\begin{proof}
  The holomorphic function $\psi^2$ is (up to a factor) \emph{the} weak Jacobi
  form of weight $-2$ and weight index $1$ (see~\cite{ez85}). Then note that
  the space of weak Jacobi forms of weight $-2$ lies inside the principal
  ideal spanned by $\psi^2$ in $J_{2*, *}$ which follows from the
  classification theorem on weak Jacobi forms in~\cite{ez85}.

  Let $X$ be any smooth projective surface with $c_1(X)^2 \neq 0$.
  By~\eqref{equ:aphigen} and the previous proposition we have
  \begin{align*}
    \phantom = & \sum_{i = 1}^\infty \ellip(A^{[[n]]})p^n
    \\
    = & \frac 1{c_1(X)^2} \left(p \frac{\diff}{\diff p}\right)^2
    \left.\frac{\diff^2}{\diff t^2}\right|_{t = 0} \ln \ellip_t(H_X(p))
    \\
    = & \psi^{-2} \sum_{i, k = 1}^\infty \sum_{m, 2l \in \IZ} i^4 k
    p^{ik}q^{mk}y^{lk} u(mi, l)
    \\
    = & \psi^{-2} \sum_{i = 1}^\infty i^4 \sum_{a|i} a^{-3} \sum_{m, 2l \in
      \IZ} u\left(\frac{mi} a, l\right) q^{am} y^{al},
  \end{align*}
  which proves the rest of the theorem.
\end{proof}

We can easily deduce G\"ottsche's and Soergel's formula on the $\chi_y$-genus
of a generalised Kummer variety (\cite{goettsche93}) from this theorem:
\begin{cor}[G\"ottsche, Soergel (\cite{goettsche93})]
  Let $n$ be a positive integer. The $\chi_y$-genus of the
  generalised Kummer variety $A^{[[n]]}$ is given by:
  \begin{gather}
    \chi_{-y}(A^{[[n]]}) = y^{n - 1} n^4 \sum_{a|n} a^{-3} \frac{y^a + y^{-a}
      - 2}{y + y^{-1} - 2}
  \end{gather}
\end{cor}

\begin{proof}
  We have
  \begin{gather*}
    \chi_{-y}(A^{[[n]]}) = y^{n - 1}
    \left.\ellip(A^{[[n]]}, \tau, z)\right|_{q = 0}
    = y^{n - 1} \left.\psi^{-2}(\tau, z)\right|_{q = 0} n^4
    \left.(\left.\psi^2\right|_{-2}
      V(n))(\tau, z)\right|_{q = 0}.
  \end{gather*}
  Now we use that $\left.\psi^2(\tau, z)\right|_{q = 0} = (y + y^{-1} - 2)$
  and the following representation of the Hecke operators (see~\cite{ez85})
  \begin{gather*}
    \left(\left.\psi^2\right|_{-2} V(n)\right)(\tau, z) = n^{-1} \sum_{a|n}
    \sum_{b = 0}^{n/a - 1} a^{-2} \psi^2\left(\frac{a \tau + b}{n/a},
      az\right),
    \\
    \intertext{i.e.}
    \begin{aligned}
      \left.\left(\left.\psi^2\right|_{-2} V(n)\right)(\tau, z)\right|_{q = 0}
      = & n^{-1} \sum_{a|n} \sum_{b = 0}^{n/a - 1} a^{-2} \lim_{\tau \to i
        \infty} \psi^2\left(\frac{a \tau + b}{n/a}, az\right)
      \\
      = & n^{-1} \sum_{a|n} \sum_{b = 0}^{n/a - 1} a^{-2} \left(y^a + y^{-a} +
        2\right)
      \\
      = & \sum_{a|n} a^{-3} \left(y^a + y^{-a} + 2\right),
    \end{aligned}
  \end{gather*}
  to conclude the formula given in the corollary.
\end{proof}

\appendix
 
\bibliographystyle{plain}
\bibliography{mybib}

\end{document}